\documentclass[11pt,reqno]{amsart}
\usepackage{amsmath,amssymb}
\usepackage[colorlinks=true,linkcolor=magenta,citecolor=blue]{hyperref}
\usepackage[margin=1in]{geometry}
\usepackage{enumitem} 
\usepackage{mathabx}
\usepackage{epsfig}

\newtheorem{dn}{Definition}[section]
\newtheorem{bdt}{Inequality}[section]
\newtheorem{dl}{Theorem}[section]
\newtheorem{md}{Proposition}[section]
\newtheorem{bd}{Lemma}[section]
\newtheorem{hq}{Corollary}[section]
\newtheorem{nx}{Remark}[section]
\newtheorem{vd}{Example}[section]
\newcommand{\R}{\mathbb{R}}

\newcommand{\e}{\varepsilon}
\newcommand{\ity}{\infty}
\newcommand{\f}{\frac}

\newcommand{\bbd}{\begin{bd}}
\newcommand{\ebd}{\end{bd}}
\newcommand{\bbdt}{\begin{bdt}}
\newcommand{\ebdt}{\end{bdt}}
\newcommand{\bdn}{\begin{dn}}
\newcommand{\edn}{\end{dn}}
\newcommand{\bhq}{\begin{hq}}
\newcommand{\ehq}{\end{hq}}
\newcommand{\bdl}{\begin{dl}}
\newcommand{\edl}{\end{dl}}
\newcommand{\bnx}{\begin{nx}}
\newcommand{\enx}{\end{nx}}
\newcommand{\bmd}{\begin{md}}
\newcommand{\emd}{\end{md}}
\newcommand{\bvd}{\begin{vd}}
\newcommand{\evd}{\end{vd}}

\title[weakly coupled system of semi-linear visco-elastic damped $\sigma$-evolution models]{A note on a weakly coupled system of semi-linear visco-elastic damped $\sigma$-evolution models with different power nonlinearities and different $\sigma$ values}
\author{Tuan Anh Dao}
\address{Tuan Anh Dao \hfill\break
$\quad$ School of Applied Mathematics and Informatics, Hanoi University of Science and Technology, No.1 Dai Co Viet road, Hanoi, Vietnam \hfill\break
Faculty for Mathematics and Computer Science, TU Bergakademie Freiberg, Pr\"{u}ferstr. 9, 09596, Freiberg, Germany}
\email{anh.daotuan@hust.edu.vn}

\begin{document}
	
\begin{abstract}
In this article, we prove the global (in time) existence of small data solutions from energy spaces basing on $L^q$ spaces, with $q \in (1,\infty)$, to the Cauchy problems for a weakly coupled system of semi-linear visco-elastic damped $\sigma$-evolution models. Here we consider different power nonlinearities and different $\sigma$ values in the comparison between two single equations. To do this, we use $(L^m \cap L^q)- L^q$ and $L^q- L^q$ estimates, i.e., by mixing additional $L^m$ regularity for the data on the basis of $L^q- L^q$ estimates for solutions, with $m \in [1,q)$, to the corresponding linear Cauchy problems. In addition, allowing loss of decay and the flexible choice of parameters $\sigma$, $m$ and $q$ bring some benefits to relax the restrictions to the admissible exponents $p$.
\end{abstract}

	\maketitle
	
	\tableofcontents
	
\section{Introduction and main results} \label{Sec.main}
Main goal of this paper is to consider the following Cauchy problems for the weakly coupled system of semi-linear visco-elastic damped $\sigma$-evolution equations:
\begin{equation}
\begin{cases}
u_{tt}+ (-\Delta)^{\sigma_1} u+ (-\Delta)^{\sigma_1} u_t=|v|^{p_1},\,\,\,  v_{tt}+ (-\Delta)^{\sigma_2} v+ (-\Delta)^{\sigma_2} v_t=|u|^{p_2}, \\
u(0,x)= u_0(x),\,\, u_t(0,x)=u_1(x),\,\, v(0,x)= v_0(x),\,\, v_t(0,x)=v_1(x), \label{pt1.1}
\end{cases}
\end{equation}
with the different values of parameters $\sigma_1,\,\sigma_2 \ge 1$ and the different power nonlinearities $p_1,\, p_2 >1$. Here the corresponding linear models with vanishing right-hand side are in the following form:
\begin{equation}
w_{tt}+ (-\Delta)^{\sigma} w+ (-\Delta)^{\sigma} w_t=0,\,\, w(0,x)= w_0(x),\,\, w_t(0,x)= w_1(x), \label{pt1.3}
\end{equation}
with $\sigma= \sigma_1$ or $\sigma= \sigma_2$. The visco-elastic damping, or so-called strong damping (see, for example, \cite{IkehataSawada,IkehataTodorova,Webb}), were studied in several recent papers. In particular, one of the most well-known models related to (\ref{pt1.3}) is those in the case $\sigma=1$ (see more \cite{Ikehata,Shibata}). Some of semi-linear visco-elastic damped wave models with power nonlinearity $|u|^p$ were investigated more detail in \cite{DabbiccoReissig} and \cite{Pizichillo}. Moreover, the authors in \cite{ReissigEbert} mentioned other interesting models to (\ref{pt1.3}) with the case $\sigma= 2$ as the visco-elastic damped plate models. Here some decay estimates of the energy and qualitative properties of solutions as well were considered. In our previous work (see \cite{DaoReissig2}) we obtained $(L^m \cap L^q)- L^q$ and $L^q- L^q$ estimates, with $q\in (1,\infty)$ and $m \in [1,q)$, for the solutions to (\ref{pt1.3}) with any $\sigma \ge 1$. To do this, there appreared two main strategies such as applying theory of modified Bessel functions combined with Fa\`{a} di Bruno's formula and the Mikhlin-H\"{o}rmander multiplier theorem to Fourier multipliers, respectively, for small frequencies and large frequencies. For this reason, the main motivation is to use these estimates in the proof of the global (in time) existence of small data energy solutions to (\ref{pt1.1}). Allowing loss of decay (see, for instance, \cite{DjaoutiReissig,DaoReissig1}) and using the fractional Gagliardo-Nirenberg inequality (see \cite{Ozawa}) come into play in the treatment of corresponding semi-linear models. Furthermore, we want to underline that how the interaction between the different values of $\sigma_1,\,\sigma_2 \ge 1$ and the flexibility of parameters $q \in (1,\ity)$, $m \in [1,q)$ affects our global (in time) existence results.\medskip

$\qquad$ \textbf{Notations}\medskip

We use the following notations throughout this paper.
\begin{itemize}[leftmargin=*]
\item We write $f\lesssim g$ when there exists a constant $C>0$ such that $f\le Cg$, and $f \approx g$ when $g\lesssim f\lesssim g$.
\item We denote $\hat{f}(t,\xi):= F_{x\rightarrow \xi}\big(f(t,x)\big)$ as the Fourier transform with respect to the space variable of a function $f(t,x)$. The spaces $H^{a,q}$ and $\dot{H}^{a,q}$, with $a \ge 0$ and $q \in (1,\ity)$, stand for Bessel and Riesz potential spaces based on $L^q$. Here $\big<D\big>^{a}$ and $|D|^{a}$ denote the pseudo-differential operators with symbols $\big<\xi\big>^{a}$ and $|\xi|^{a}$, respectively.
\item We fix the following constants with $q \in (1,\ity)$, $m \in [1,q)$ and $n \ge 1$:
\begin{align*}
&\alpha:= \Big(\frac{1}{m}-\frac{1}{q}\Big)\Big(2+\Big[\frac{n}{2}\Big]+ n\Big(\frac{1}{\sigma_2}- \frac{1}{\sigma_1}\Big)\Big), \\
&\beta:= \Big(\frac{1}{m}-\frac{1}{q}\Big)\Big(2+\Big[\frac{n}{2}\Big]+ n\Big(\frac{1}{\sigma_1}- \frac{1}{\sigma_2}\Big)\Big), \\
&\gamma:= \Big(1+\frac{1}{q}-\frac{1}{m}\Big)\Big(2+\Big[\frac{n}{2}\Big]\Big).
\end{align*}
Then, we denote $\kappa_1:= \frac{1}{2}(1+\gamma+\alpha)$ and $\kappa_2:= \frac{1}{2}(1+\gamma+\beta)$.
\item Finally, we introduce the spaces
$\mathcal{A}^{j}_{m,q}:= \big(L^m \cap H^{2\sigma_j,q}\big) \times \big(L^m \cap L^q\big)$ with the norm
$$\|(u_0,u_1)\|_{\mathcal{A}^{j}_{m,q}}:=\|u_0\|_{L^m}+ \|u_0\|_{H^{2\sigma_j,q}}+ \|u_1\|_{L^m}+ \|u_1\|_{L^q}, $$
where $q \in (1,\ity)$, $m \in [1,q)$ and $j=1,\,2$.
\end{itemize}

$\qquad$ \textbf{Main results}\medskip

Let us state the main results which will be proved in the present paper.
\bdl[\textbf{Loss of decay}] \label{dl1.1}
Let $q \in (1,\ity)$ be a fixed constant and $m\in [1,q)$. Let assume $\sigma_1 \ge \sigma_2$ and $n> \sigma_1$. We assume the conditions
\begin{align}
&\frac{q}{m} \le p_1,\, p_2 < \ity & & &\quad \text{ if }&\, n \le 2q\sigma_2, \label{GN11A1} \\
&\frac{q}{m} \le p_1 \le \frac{n}{n- 2q\sigma_2}, &\quad &\frac{q}{m} \le p_2  < \ity &\quad \text{ if }&\, 2q\sigma_2 < n \le 2q\sigma_1, \label{GN11A2} \\
&\frac{q}{m} \le p_1 \le \frac{n}{n- 2q\sigma_2}, &\quad &\frac{q}{m} \le p_2  \le \frac{n}{n- 2q\sigma_1} &\quad \text{ if }&\, 2q\sigma_1 < n \le \frac{2q^2\sigma_2}{q-m}. \label{GN11A3}
\end{align}
Moreover, we suppose the following conditions:
\begin{equation} \label{exponent11A}
m \Big(\frac{1+ p_2 +p_2(1+p_1)\kappa_1}{\frac{p_2-1}{\sigma_1}+ p_2\frac{p_1- 1}{\sigma_2}}\Big) < \frac{n}{2}, \text{ and }p_1 \le 1+\frac{2m\sigma_2(1+\kappa_1)}{n- 2m\sigma_2\kappa_1} < p_2.
\end{equation}
Then, there exists a constant $\e>0$ such that for any small data
\begin{equation*}
\big((u_0,u_1),\, (v_0,v_1) \big) \in \mathcal{A}^{1}_{m,q} \times \mathcal{A}^{2}_{m,q} \text{ satisfying the assumption } \|(u_0,u_1)\|_{\mathcal{A}^{1}_{m,q}}+ \|(v_0,v_1)\|_{\mathcal{A}^{2}_{m,q}} \le \e,
\end{equation*}
we have a uniquely determined global (in time) small data energy solution
$$ (u,v) \in \Big(C([0,\ity),H^{2\sigma_1,q})\cap C^1([0,\ity),L^q\Big) \times \Big(C([0,\ity),H^{2\sigma_2,q})\cap C^1([0,\ity),L^q\Big) $$
to (\ref{pt1.1}). The following estimates hold:
\begin{align}
\|u(t,\cdot)\|_{L^q}& \lesssim (1+t)^{-\frac{n}{2\sigma_1}(1-\frac{1}{r})+ \e(p_1,\sigma_2)+ \frac{\kappa_1}{2}} \big(\|(u_0,u_1)\|_{\mathcal{A}^{1}_{m,q}}+ \|(v_0,v_1)\|_{\mathcal{A}^{2}_{m,q}}\big), \label{decayrate11A1} \\
\big\|\big(|D|^{\sigma_1} u(t,\cdot),u_t(t,\cdot)\big)\big\|_{L^q}& \lesssim (1+t)^{-\frac{n}{2\sigma_1}(1-\frac{1}{r})- \frac{1}{2}+ \e(p_1,\sigma_2)+ \frac{\kappa_1}{2}} \big(\|(u_0,u_1)\|_{\mathcal{A}^{1}_{m,q}}+ \|(v_0,v_1)\|_{\mathcal{A}^{2}_{m,q}}\big), \label{decayrate11A2} \\
\big\||D|^{2\sigma_1} u(t,\cdot)\big\|_{L^q}& \lesssim (1+t)^{-\frac{n}{2\sigma_1}(1-\frac{1}{r})-1+ \e(p_1,\sigma_2) + \frac{\kappa_1}{2}} \big(\|(u_0,u_1)\|_{\mathcal{A}^{1}_{m,q}}+ \|(v_0,v_1)\|_{\mathcal{A}^{2}_{m,q}}\big), \label{decayrate11A3} \\
\|v(t,\cdot)\|_{L^q}& \lesssim (1+t)^{-\frac{n}{2\sigma_2}(1-\frac{1}{r})+ \frac{\kappa_1}{2}} \big(\|(u_0,u_1)\|_{\mathcal{A}^{1}_{m,q}}+ \|(v_0,v_1)\|_{\mathcal{A}^{2}_{m,q}}\big), \label{decayrate11A4} \\
\big\|\big(|D|^{\sigma_2} v(t,\cdot),v_t(t,\cdot)\big)\big\|_{L^q}& \lesssim (1+t)^{-\frac{n}{2\sigma_2}(1-\frac{1}{r})- \frac{1}{2}+ \frac{\kappa_1}{2}} \big(\|(u_0,u_1)\|_{\mathcal{A}^{1}_{m,q}}+ \|(v_0,v_1)\|_{\mathcal{A}^{2}_{m,q}}\big), \label{decayrate11A5} \\
\big\||D|^{2\sigma_2} v(t,\cdot)\big\|_{L^q}& \lesssim (1+t)^{-\frac{n}{2\sigma_2}(1-\frac{1}{r})-1+ \frac{\kappa_1}{2}} \big(\|(u_0,u_1)\|_{\mathcal{A}^{1}_{m,q}}+ \|(v_0,v_1)\|_{\mathcal{A}^{2}_{m,q}}\big), \label{decayrate11A6}
\end{align}
where $\e(p_1,\sigma_2):= 1- \frac{n}{2m\sigma_2}(p_1- 1)+ p_1\kappa_1$.
\edl

\bdl[\textbf{Loss of decay}] \label{dl1.2}
Let $q \in (1,\ity)$ be a fixed constant and $m\in [1,q)$. Let assume $\sigma_2 \ge \sigma_1$ and $n> \sigma_2$. We assume the conditions
\begin{align}
&\frac{q}{m} \le p_1,\, p_2 < \ity & & &\quad \text{ if }&\, n \le 2q\sigma_1, \label{GN12A1} \\
&\frac{q}{m} \le p_1 < \ity, &\quad &\frac{q}{m} \le p_2  \le \frac{n}{n- 2q\sigma_1} &\quad \text{ if }&\, 2q\sigma_1 < n \le 2q\sigma_2, \label{GN12A2} \\
&\frac{q}{m} \le p_1 \le \frac{n}{n- 2q\sigma_2}, &\quad &\frac{q}{m} \le p_2  \le \frac{n}{n- 2q\sigma_1} &\quad \text{ if }&\, 2q\sigma_2 < n \le \frac{2q^2\sigma_1}{q-m}. \label{GN12A3}
\end{align}
Moreover, we suppose the following conditions:
\begin{equation} \label{exponent12A}
m \Big(\frac{1+ p_1 +p_1(1+p_2)\kappa_2}{\frac{p_1-1}{\sigma_2}+ p_1\frac{p_2- 1}{\sigma_1}}\Big) < \frac{n}{2}, \text{ and }p_2 \le 1+\frac{2m\sigma_1(1+\kappa_2)}{n- 2m\sigma_1\kappa_2} < p_1.
\end{equation}
Then, there exists a constant $\e>0$ such that for any small data
\begin{equation*}
\big((u_0,u_1),\, (v_0,v_1) \big) \in \mathcal{A}^{1}_{m,q} \times \mathcal{A}^{2}_{m,q} \text{ satisfying the assumption } \|(u_0,u_1)\|_{\mathcal{A}^{1}_{m,q}}+ \|(v_0,v_1)\|_{\mathcal{A}^{2}_{m,q}} \le \e,
\end{equation*}
we have a uniquely determined global (in time) small data energy solution
$$ (u,v) \in \Big(C([0,\ity),H^{2\sigma_1,q})\cap C^1([0,\ity),L^q\Big) \times \Big(C([0,\ity),H^{2\sigma_2,q})\cap C^1([0,\ity),L^q\Big) $$
to (\ref{pt1.1}). The following estimates hold:
\begin{align}
\|u(t,\cdot)\|_{L^q}& \lesssim (1+t)^{-\frac{n}{2\sigma_1}(1-\frac{1}{r})+ \frac{\kappa_2}{2}} \big(\|(u_0,u_1)\|_{\mathcal{A}^{1}_{m,q}}+ \|(v_0,v_1)\|_{\mathcal{A}^{2}_{m,q}}\big), \label{decayrate12A1} \\
\big\|\big(|D|^{\sigma_1} u(t,\cdot), u_t(t,\cdot)\big)\big\|_{L^q}& \lesssim (1+t)^{-\frac{n}{2\sigma_1}(1-\frac{1}{r})- \frac{1}{2}+ \frac{\kappa_2}{2}} \big(\|(u_0,u_1)\|_{\mathcal{A}^{1}_{m,q}}+ \|(v_0,v_1)\|_{\mathcal{A}^{2}_{m,q}}\big), \label{decayrate12A2} \\
\big\||D|^{2\sigma_1} u(t,\cdot)\big\|_{L^q}& \lesssim (1+t)^{-\frac{n}{2\sigma_1}(1-\frac{1}{r})-1+ \frac{\kappa_2}{2}} \big(\|(u_0,u_1)\|_{\mathcal{A}^{1}_{m,q}}+ \|(v_0,v_1)\|_{\mathcal{A}^{2}_{m,q}}\big), \label{decayrate12A3} \\
\|v(t,\cdot)\|_{L^q}& \lesssim (1+t)^{-\frac{n}{2\sigma_2}(1-\frac{1}{r})+ \e(p_2,\sigma_1)+ \frac{\kappa_2}{2}} \big(\|(u_0,u_1)\|_{\mathcal{A}^{1}_{m,q}}+ \|(v_0,v_1)\|_{\mathcal{A}^{2}_{m,q}}\big), \label{decayrate12A4} \\
\big\|\big(|D|^{\sigma_2} v(t,\cdot), v_t(t,\cdot)\big)\big\|_{L^q}& \lesssim (1+t)^{-\frac{n}{2\sigma_2}(1-\frac{1}{r})- \frac{1}{2}+ \e(p_2,\sigma_1)+ \frac{\kappa_2}{2}} \big(\|(u_0,u_1)\|_{\mathcal{A}^{1}_{m,q}}+ \|(v_0,v_1)\|_{\mathcal{A}^{2}_{m,q}}\big), \label{decayrate12A5} \\
\big\||D|^{2\sigma_2} v(t,\cdot)\big\|_{L^q}& \lesssim (1+t)^{-\frac{n}{2\sigma_2}(1-\frac{1}{r})-1+ \e(p_2,\sigma_1)+ \frac{\kappa_2}{2}} \big(\|(u_0,u_1)\|_{\mathcal{A}^{1}_{m,q}}+ \|(v_0,v_1)\|_{\mathcal{A}^{2}_{m,q}}\big), \label{decayrate12A6}
\end{align}
where $\e(p_2,\sigma_1):= 1- \frac{n}{2m\sigma_1}(p_2- 1)+ p_2\kappa_2$.
\edl

\bdl[\textbf{No loss of decay}] \label{dl1.1B}
Under the analogous assumptions of Theorem \ref{dl1.1}, if condition (\ref{exponent11A}) is replaced by
\begin{equation} \label{exponent11B}
\min\{p_1,\,p_2\}> 1+\frac{2m\sigma_2(1+\kappa_1)}{n- 2m\sigma_2\kappa_1},
\end{equation}
then we have the same conclusions of Theorem \ref{dl1.1}. But the estimates (\ref{decayrate11A1})-(\ref{decayrate11A6}) are modified in the following way:
\begin{align}
\|u(t,\cdot)\|_{L^q}& \lesssim (1+t)^{-\frac{n}{2\sigma_1}(1-\frac{1}{r})+ \frac{\gamma+1}{2}} \big(\|(u_0,u_1)\|_{\mathcal{A}^{1}_{m,q}}+ \|(v_0,v_1)\|_{\mathcal{A}^{2}_{m,q}}\big), \label{decayrate11B1} \\
\big\|\big(|D|^{\sigma_1} u(t,\cdot),u_t(t,\cdot)\big)\big\|_{L^q}& \lesssim (1+t)^{-\frac{n}{2\sigma_1}(1-\frac{1}{r})+ \frac{\gamma}{2}} \big(\|(u_0,u_1)\|_{\mathcal{A}^{1}_{m,q}}+ \|(v_0,v_1)\|_{\mathcal{A}^{2}_{m,q}}\big), \label{decayrate11B2} \\
\big\||D|^{2\sigma_1} u(t,\cdot)\big\|_{L^q}& \lesssim (1+t)^{-\frac{n}{2\sigma_1}(1-\frac{1}{r})+ \frac{\gamma-1}{2}} \big(\|(u_0,u_1)\|_{\mathcal{A}^{1}_{m,q}}+ \|(v_0,v_1)\|_{\mathcal{A}^{2}_{m,q}}\big), \label{decayrate11B3} \\
\|v(t,\cdot)\|_{L^q}& \lesssim (1+t)^{-\frac{n}{2\sigma_2}(1-\frac{1}{r})+ \frac{\gamma+1}{2}} \big(\|(u_0,u_1)\|_{\mathcal{A}^{1}_{m,q}}+ \|(v_0,v_1)\|_{\mathcal{A}^{2}_{m,q}}\big), \label{decayrate11B4} \\
\big\|\big(|D|^{\sigma_2} v(t,\cdot),v_t(t,\cdot)\big)\big\|_{L^q}& \lesssim (1+t)^{-\frac{n}{2\sigma_2}(1-\frac{1}{r})+ \frac{\gamma}{2}} \big(\|(u_0,u_1)\|_{\mathcal{A}^{1}_{m,q}}+ \|(v_0,v_1)\|_{\mathcal{A}^{2}_{m,q}}\big), \label{decayrate11B5} \\
\big\||D|^{2\sigma_2} v(t,\cdot)\big\|_{L^q}& \lesssim (1+t)^{-\frac{n}{2\sigma_2}(1-\frac{1}{r})+ \frac{\gamma-1}{2}} \big(\|(u_0,u_1)\|_{\mathcal{A}^{1}_{m,q}}+ \|(v_0,v_1)\|_{\mathcal{A}^{2}_{m,q}}\big). \label{decayrate11B6}
\end{align}
\edl

\bdl[\textbf{No loss of decay}] \label{dl1.2B}
Under the analogous assumptions of Theorem \ref{dl1.2}, if condition (\ref{exponent12A}) is replaced by
\begin{equation} \label{exponent12B}
\min\{p_1,\,p_2\}> 1+\frac{2m\sigma_1(1+\kappa_2)}{n- 2m\sigma_1\kappa_2},
\end{equation}
then we have the same conclusions of Theorem \ref{dl1.2}. Moreover, the estimates (\ref{decayrate11B1})-(\ref{decayrate11B6}) hold instead of (\ref{decayrate12A1})-(\ref{decayrate12A6}).
\edl

\begin{nx}
\fontshape{n}
\selectfont
Let us observe the interplay between the parameters $\sigma_j$ and $p_j$ with $j=1,2$. In Theorems \ref{dl1.1} to \ref{dl1.2B}, we can see that the different choice of $\sigma_1$ and $\sigma_2$ influences our admissible exponents $p_1$ and $p_2$. Moreover, if we want to consider the special case of $\sigma_1= \sigma_2= \sigma$ then it is clear that the constants $\kappa_1= \kappa_2= \kappa$. For this reason, the conditions from (\ref{GN11A1}) to (\ref{GN11A3}) are similar to those from (\ref{GN12A1}) to (\ref{GN12A3}). We also re-write both the assumptions (\ref{exponent11A}) and (\ref{exponent12A}) in the following common form:
$$ m \sigma \Big(\frac{1+ \max\{p_1,\,p_2\}+ (p_1p_2+ \max\{p_1,\,p_2\})\kappa}{p_1p_2- 1}\Big) < \frac{n}{2}, $$
and
$$ \min\{p_1,\,p_2\} \le 1+\frac{2m\sigma(1+\kappa)}{n- 2m\sigma\kappa} < \max\{p_1,\,p_2\}. $$
Finally, in this case we may conclude that all the statements in Theorems \ref{dl1.1B} and \ref{dl1.2B} coincide.
\end{nx}

\begin{nx}
\fontshape{n}
\selectfont
Here we want to underline that due to the conditions (\ref{exponent11A}) and (\ref{exponent12A}), both $\e(p_1,\sigma_2)$ and $\e(p_2,\sigma_1)$ in Theorems \ref{dl1.1} and \ref{dl1.2} are non-negative.
\end{nx}

\begin{nx}
\fontshape{n}
\selectfont
Let us explain our main strategies appearing in Theorems \ref{dl1.1} to \ref{dl1.2B} such as loss of decay and no loss of decay. In particular, the first phenomenon means that the decay rates for the solutions to (\ref{pt1.1}) in Theorems \ref{dl1.1} and \ref{dl1.2} are worse than those for solutions to (\ref{pt1.3}). Meanwhile, the second one shows that in comparison with the corresponding linear models these decay rates in Theorems \ref{dl1.1B} and \ref{dl1.2B} are exactly the same. It is clear that allowing loss of decay comes into play to relax the restrictions to the admissible exponents $p_1$ and $p_2$. Indeed, here we can see that we allow one exponent $p_1$ in (\ref{exponent11A}) or $p_2$ in (\ref{exponent12A}) below the exponent $1+\frac{2m\sigma_2(1+\kappa_1)}{n- 2m\sigma_2\kappa_1}$ or $1+\frac{2m\sigma_1(1+\kappa_2)}{n- 2m\sigma_1\kappa_2}$, respectively. On contrary, in (\ref{exponent11B}) and (\ref{exponent12B}) we need to guarantee both exponents $p_1$ and $p_2$ above these exponents.
\end{nx}

\textbf{This article is organized as follows}: In Section \ref{Sec2}, we collect $(L^m \cap L^q)- L^q$ and $L^q- L^q$ estimates for solutions to (\ref{pt1.3}), with $q\in (1,\infty)$ and $m \in [1,q)$, from our previous work \cite{DaoReissig2} and prove a sharper result as well. We give the proofs of our global (in time) existence results to (\ref{pt1.1}) in Section \ref{Semi-linear estimates}.

\section{Preliminaries} \label{Sec2}
In this section, main goal is to derive $(L^m \cap L^q)- L^q$ and $L^q- L^q$ estimates for the solution and some its derivatives to (\ref{pt1.3}). First, using partial Fourier transformation to (\ref{pt1.3}) we have the following Cauchy problem:
\begin{equation}
\hat{w}_{tt}+ |\xi|^{2\sigma} \hat{w}_t+ |\xi|^{2\sigma} \hat{w}=0,\,\, \hat{w}(0,\xi)= \hat{w}_0(\xi),\,\, \hat{w}_t(0,\xi)= \hat{w}_1(\xi). \label{pt3.1}
\end{equation}
The characteristic roots are
$$ \lambda_{1,2}=\lambda_{1,2}(\xi)= \f{1}{2}\Big(-|\xi|^{2\sigma}\pm \sqrt{|\xi|^{4\sigma}-4|\xi|^{2\sigma}}\Big). $$
We write the solutions to (\ref{pt3.1}) in the following form:
\begin{equation}
\hat{w}(t,\xi)= \frac{\lambda_1 e^{\lambda_2 t}-\lambda_2 e^{\lambda_1 t}}{\lambda_1- \lambda_2}\hat{w}_0(\xi)+ \frac{e^{\lambda_1 t}-e^{\lambda_2 t}}{\lambda_1- \lambda_2}\hat{w}_1(\xi)=: \hat{K}_{0,\sigma}(t,\xi)\hat{w}_0(\xi)+\hat{K}_{1,\sigma}(t,\xi)\hat{w}_1(\xi), \label{pt3.2}
\end{equation}
where we assume $\lambda_{1}\neq \lambda_{2}$. Taking account of the cases of small and large frequencies separately, we get the following asymptotic behavior of the characteristic roots:
\begin{align}
&1.\,\, \lambda_{1,2} \sim -|\xi|^{2\sigma}\pm i|\xi|^\sigma,\,\, \lambda_1-\lambda_2 \sim i|\xi|^\sigma \text{ for small } |\xi|, \label{pt3.3} \\
&2.\,\, \lambda_1\sim -1,\,\, \lambda_2\sim -|\xi|^{2\sigma},\,\, \lambda_1-\lambda_2 \sim |\xi|^{2\sigma} \text{ for large } |\xi|. \label{pt3.4}
\end{align}
We now re-write the solution to (\ref{pt1.3}) which are decomposed into two parts localized separately to low and high frequencies in the following form:
$$ w(t,x)= w_\chi(t,x)+ w_{1-\chi}(t,x), $$
where
$$w_\chi(t,x)=F^{-1}_{\xi\rightarrow x}\big(\chi(|\xi|)\hat{w}(t,\xi)\big) \text{ and } w_{1-\chi}(t,x)=F^{-1}_{\xi\rightarrow x}\big(\big(1-\chi(|\xi|)\big)\hat{w}(t,\xi)\big). $$
Here $\chi(|\xi|)$ is a smooth cut-off function equal to $1$ for small $|\xi|$ and vanishing for large $|\xi|$. Recalling the statements from Proposition $3.9$ in \cite{DaoReissig2}, we obtained the following result.

\bmd \label{md2.1}
Let $\sigma \ge 1$ in (\ref{pt1.3}), $q\in (1,\ity)$ and $m\in [1,q)$. Then the Sobolev solutions to (\ref{pt1.3}) satisfy the following $(L^m \cap L^q)-L^q$ estimates:
\begin{align*}
\big\||D|^a w(t,\cdot)\big\|_{L^q}& \lesssim  (1+t)^{\frac{1}{2}(2+[\frac{n}{2}])\frac{1}{r} -\frac{n}{2\sigma}(1-\frac{1}{r})-\frac{a}{2\sigma}} \|w_0\|_{L^m \cap H^{a,q}} \\
& \qquad \quad + (1+t)^{1+\frac{1}{2}(1+[\frac{n}{2}])\frac{1}{r} -\frac{n}{2\sigma}(1-\frac{1}{r})-\frac{a}{2\sigma}}\|w_1\|_{L^m \cap H^{[a-2\sigma]^+,q}}, \\
\big\||D|^a w_t(t,\cdot)\big\|_{L^q}& \lesssim  (1+t)^{\frac{1}{2}(1+[\frac{n}{2}])\frac{1}{r} -\frac{n}{2\sigma}(1-\frac{1}{r})-\frac{a}{2\sigma}} \|w_0\|_{L^m \cap H^{a,q}} \\
& \qquad \quad + (1+t)^{\frac{1}{2}(2+[\frac{n}{2}])\frac{1}{r} -\frac{n}{2\sigma}(1-\frac{1}{r})-\frac{a}{2\sigma}}\|w_1\|_{L^m \cap H^{a,q}}.
\end{align*}
Moreover, the following $L^q-L^q$ estimates hold:
\begin{align*}
\big\||D|^a w(t,\cdot)\big\|_{L^q}& \lesssim (1+t)^{\frac{1}{2}(2+[\frac{n}{2}])-\frac{a}{2\sigma}} \|w_0\|_{H^{a,q}}+ (1+t)^{1+\frac{1}{2}(1+[\frac{n}{2}])-\frac{a}{2\sigma}} \|w_1\|_{H^{[a-2\sigma]^+,q}}, \\
\big\||D|^a w_t(t,\cdot)\big\|_{L^q}& \lesssim (1+t)^{\frac{1}{2}(1+[\frac{n}{2}])-\frac{a}{2\sigma}} \|w_0\|_{H^{a,q}}+ (1+t)^{\frac{1}{2}(2+[\frac{n}{2}])-\frac{a}{2\sigma}} \|w_1\|_{H^{a,q}}.
\end{align*}
Here $1+ \frac{1}{q}= \frac{1}{r}+\frac{1}{m}$, $a$ is a non-negative number and for all the dimension $n\ge 1$.
\emd

\begin{nx}
\fontshape{n}
\selectfont
In Proposition \ref{md2.1} we derived estimates for the solution and some its derivatives to (\ref{pt1.3}) for any space dimensions $n \ge 1$. In addition, we want to underline that under a constraint condition to space dimensions $n> \sigma$ we may prove the following sharper result.
\end{nx}

\bmd \label{md2.2}
Let $\sigma \ge 1$ in (\ref{pt1.3}), $q\in (1,\ity)$ and $m\in [1,q)$. Then the Sobolev solutions to (\ref{pt1.3}) satisfy the following $(L^m \cap L^q)-L^q$ estimates:
\begin{align*}
\big\||D|^a w(t,\cdot)\big\|_{L^q}& \lesssim  (1+t)^{\frac{1}{2}(2+[\frac{n}{2}])\frac{1}{r} -\frac{n}{2\sigma}(1-\frac{1}{r})-\frac{a}{2\sigma}} \|w_0\|_{L^m \cap H^{a,q}} \\
& \qquad \quad + (1+t)^{\frac{1}{2}+\frac{1}{2}(2+[\frac{n}{2}])\frac{1}{r} -\frac{n}{2\sigma}(1-\frac{1}{r})-\frac{a}{2\sigma}}\|w_1\|_{L^m \cap H^{[a-2\sigma]^+,q}}, \\
\big\||D|^a w_t(t,\cdot)\big\|_{L^q}& \lesssim  (1+t)^{\frac{1}{2}(2+[\frac{n}{2}])\frac{1}{r} -\frac{n}{2\sigma}(1-\frac{1}{r})-\frac{1}{2}-\frac{a}{2\sigma}} \|w_0\|_{L^m \cap H^{a,q}} \\
& \qquad \quad + (1+t)^{\frac{1}{2}(2+[\frac{n}{2}])\frac{1}{r} -\frac{n}{2\sigma}(1-\frac{1}{r})-\frac{a}{2\sigma}}\|w_1\|_{L^m \cap H^{a,q}}.
\end{align*}
Moreover, the following $L^q-L^q$ estimates hold:
\begin{align*}
\big\||D|^a w(t,\cdot)\big\|_{L^q}& \lesssim (1+t)^{\frac{1}{2}(2+[\frac{n}{2}])-\frac{a}{2\sigma}} \|w_0\|_{H^{a,q}}+ (1+t)^{\frac{1}{2}(3+[\frac{n}{2}])-\frac{a}{2\sigma}} \|w_1\|_{H^{[a-2\sigma]^+,q}}, \\
\big\||D|^a w_t(t,\cdot)\big\|_{L^q}& \lesssim (1+t)^{\frac{1}{2}(1+[\frac{n}{2}])- \frac{a}{2\sigma}} \|w_0\|_{H^{a,q}}+ (1+t)^{\frac{1}{2}(2+[\frac{n}{2}])-\frac{a}{2\sigma}} \|w_1\|_{H^{a,q}}.
\end{align*}
Here $1+ \frac{1}{q}= \frac{1}{r}+\frac{1}{m}$, $a$ is a non-negative number and for all the dimension $n> \sigma$.
\emd

In order to prove Proposition \ref{md2.2}, we shall show the following auxiliary estimates.
\bbd \label{bd2.1}
The following estimates hold in $\R^n$ for any $n> \sigma$:
\begin{align*}
\big\|F^{-1}\big(\chi(|\xi|)|\xi|^a \hat{K}_{0,\sigma}(t,\xi)\big)(t,\cdot)\big\|_{L^r}&\lesssim
\begin{cases}
1 \text{ if } t\in (0,1], &\\
t^{\frac{1}{2}(2+[\frac{n}{2}])\frac{1}{r} -\frac{n}{2\sigma}(1-\frac{1}{r})-\frac{a}{2\sigma}} \text{ if } t\in[1,\ity),&
\end{cases} \\
\big\|F^{-1}\big(\chi(|\xi|)|\xi|^a \hat{K}_{1,\sigma}(t,\xi)\big)(t,\cdot)\big\|_{L^r}&\lesssim
\begin{cases}
1 \text{ if } t\in (0,1], &\\
t^{\frac{1}{2}+\frac{1}{2}(2+[\frac{n}{2}])\frac{1}{r} -\frac{n}{2\sigma}(1-\frac{1}{r})-\frac{a}{2\sigma}} \text{ if } t\in[1,\ity), &
\end{cases}
\end{align*}
for all $r \in [1,\ity]$ and any non-negative number $a$.
\ebd

\begin{proof}
The proof of the first statement is completed from Proposition $3.7$ in \cite{DaoReissig2}. For this reason, we will prove the second one only. Thanks to the asymptotic behavior of the characteristic roots in (\ref{pt3.3}), we arrive at immediately the following estimate for small frequencies:
$$ \big|\hat{K}_{1,\sigma}(t,\xi)\big|= \frac{\big|e^{\lambda_1 t}-e^{\lambda_2 t}\big|}{|\lambda_1- \lambda_2|} \lesssim |\xi|^{-\sigma} e^{-c|\xi|^{2\sigma}t}, $$
where $c$ is a suitable positive constant. We can see that it holds for small frequencies
\begin{equation*}
\int_{\R^n} |\xi|^{b_1} e^{-c|\xi|^{b_2} t}d\xi \lesssim (1+t)^{-\frac{n+b_1}{b_2}},
\end{equation*}
for any $n \ge 1$, $b_1 \in \R$ satisfying $n+b_1 >0$ and for all positive numbers $c,\,b_2 >0$. Hence, we derive
\begin{equation} \label{l2.1.1}
\big\|F^{-1}\big(\chi(|\xi|)|\xi|^a \hat{K}_{1,\sigma}(t,\xi)\big)(t,\cdot)\big\|_{L^\ity} \lesssim \int_{\R^n} \chi(|\xi|) e^{-|\xi|^{2\sigma}t} |\xi|^{a-\sigma} d\xi \lesssim (1+t)^{-\frac{n+a-\sigma}{2\sigma}}.
\end{equation}
From Proposition $3.1$ and Remark $3.3$ in \cite{DaoReissig2}, we get
\begin{equation} \label{l2.1.2}
\big\| F^{-1}\big(\chi(|\xi|)|\xi|^a \hat{K}_{1,\sigma}(t,\xi)\big)(t,\cdot)\big\|_{L^1}\lesssim \begin{cases}
t \text{ for } t\in (0,1], &\\
t^{1+\frac{1}{2}(1+[\frac{n}{2}])-\frac{a}{2\sigma}} \text{ for } t\in[1,\ity). &
\end{cases}
\end{equation}
By interpolation argument, from (\ref{l2.1.1}) and (\ref{l2.1.2}) we may conclude the second statement that we wanted to prove. 
\end{proof}

\begin{proof}[Proof of Proposition \ref{md2.2}]
In order to obtain the $(L^m \cap L^q)- L^q$ estimates, we combine the statements from Lemma \ref{bd2.1} with those from Corollary $3.4$ in \cite{DaoReissig2}. In particular, we can estimate the $L^q$ norm of the small frequency part of the solutions by the $L^m$ norm of the data, whereas its high-frequency part is controlled by using the $L^q-L^q$ estimates. Finally, applying Young's convolution inequality we may conclude all the statements from Proposition \ref{md2.2}. This completes our proof.
\end{proof}

\noindent From the statements in Proposition \ref{md2.2}, replacing $\sigma= \sigma_j$ with $j=1,2$ and recalling abbreviation $\gamma= (2+[\frac{n}{2}])\frac{1}{r}$ we have the following result.
\bhq \label{hq2.1}
Let $\sigma= \sigma_j \ge 1$ with $j=1,\,2$ in (\ref{pt1.3}). Let $q\in (1,\ity)$ and $m\in [1,q)$. Then the Sobolev solutions to (\ref{pt1.3}) satisfy the following $(L^m \cap L^q)-L^q$ estimates:
\begin{align*}
\big\||D|^a w(t,\cdot)\big\|_{L^q}& \lesssim  (1+t)^{\frac{\gamma}{2} -\frac{n}{2\sigma_j}(1-\frac{1}{r})-\frac{a}{2\sigma_j}} \|w_0\|_{L^m \cap H^{a,q}} \\
& \qquad \quad + (1+t)^{\frac{\gamma+1}{2} -\frac{n}{2\sigma_j}(1-\frac{1}{r})-\frac{a}{2\sigma_j}}\|w_1\|_{L^m \cap H^{[a-2\sigma_j]^+,q}}, \\
\big\||D|^a w_t(t,\cdot)\big\|_{L^q}& \lesssim  (1+t)^{\frac{\gamma-1}{2} -\frac{n}{2\sigma_j}(1-\frac{1}{r})-\frac{a}{2\sigma_j}} \|w_0\|_{L^m \cap H^{a,q}} \\
& \qquad \quad + (1+t)^{\frac{\gamma}{2} -\frac{n}{2\sigma_j}(1-\frac{1}{r})-\frac{a}{2\sigma_j}}\|w_1\|_{L^m \cap H^{a,q}}.
\end{align*}
Moreover, the following $L^q-L^q$ estimates hold:
\begin{align*}
\big\||D|^a w(t,\cdot)\big\|_{L^q}& \lesssim (1+t)^{\frac{1}{2}(2+[\frac{n}{2}])-\frac{a}{2\sigma_j}} \|w_0\|_{H^{a,q}}+ (1+t)^{\frac{1}{2}(3+[\frac{n}{2}])-\frac{a}{2\sigma_j}} \|w_1\|_{H^{[a-2\sigma_j]^+,q}}, \\
\big\||D|^a w_t(t,\cdot)\big\|_{L^q}& \lesssim (1+t)^{\frac{1}{2}(1+[\frac{n}{2}])- \frac{a}{2\sigma_j}} \|w_0\|_{H^{a,q}}+ (1+t)^{\frac{1}{2}(2+[\frac{n}{2}])-\frac{a}{2\sigma_j}} \|w_1\|_{H^{a,q}}.
\end{align*}
Here $1+ \frac{1}{q}= \frac{1}{r}+\frac{1}{m}$, $a$ is a non-negative number and for all the dimension $n> \max\{\sigma_1,\, \sigma_2\}$.
\ehq

\begin{nx}
\fontshape{n}
\selectfont
The $(L^m \cap L^q)- L^q$ and $L^q- L^q$ estimates for solutions to (\ref{pt1.3}), with $q\in (1,\infty)$ and $m \in [1,q)$ come into play in treatment of the weakly coupled system of corresponding semi-linear models (\ref{pt1.1}) in next section. It is clear that the decay estimates for solution and some its derivatives to (\ref{pt1.3}) from Proposition \ref{md2.2} are better than those from Proposition \ref{md2.1}. Hence, it is reasonable to apply the statements from Corollary \ref{hq2.1} in the steps of the proofs to our global (in time) existence results. 
\end{nx}

\section{Proofs of the global (in time) existence results} \label{Semi-linear estimates}

\subsection{Proof of Theorem \ref{dl1.1}}
Recalling the fundamental solutions $K_{0,\sigma}$ and $K_{1,\sigma}$ defined in Section \ref{Sec2} we write the solutions of the corresponding linear Cauchy problems with vanishing right-hand sides to (\ref{pt1.1}) as follows:
$$\begin{cases}
u^{ln}(t,x)=K_{0,\sigma_1}(t,x) \ast_{x} u_0(x)+ K_{1,\sigma_1}(t,x) \ast_{x} u_1(x), \\ 
v^{ln}(t,x)=K_{0,\sigma_2}(t,x) \ast_{x} v_0(x)+ K_{1,\sigma_2}(t,x) \ast_{x} v_1(x).
\end{cases}$$
Using Duhamel's principle we get the formal implicit representation of the solutions to (\ref{pt1.1}) in the following form:
$$\begin{cases}
u(t,x)= u^{ln}(t,x) + \int_0^t K_{1,\sigma_1}(t-\tau,x) \ast_x |v(t,x)|^{p_1} d\tau=: u^{ln}(t,x)+ u^{nl}(t,x), \\ 
v(t,x)= v^{ln}(t,x) + \int_0^t K_{1,\sigma_2}(t-\tau,x) \ast_x |u(t,x)|^{p_2} d\tau=: v^{ln}(t,x)+ v^{nl}(t,x).
\end{cases}$$
First, we choose the data spaces $(u_0,u_1) \in \mathcal{A}^{1}_{m,q}$ and $(v_0,v_1) \in \mathcal{A}^{2}_{m,q}$. We introduce the family $\{X(t)\}_{t>0}$ of solution spaces $X(t)$ with the norm
\begin{align*}
\|(u,v)\|_{X(t)}:= \sup_{0\le \tau \le t} \Big( f_{1}(\tau)^{-1}\|u(\tau,\cdot)\|_{L^q} &+ f_{1,\sigma_1}(\tau)^{-1}\big\||D|^{\sigma_1} u(\tau,\cdot)\big\|_{L^q}+ f_{1,2\sigma_1}(\tau)^{-1}\big\||D|^{2\sigma_1} u(\tau,\cdot)\big\|_{L^q}\\
&+ f_{2}(\tau)^{-1}\|u_t(\tau,\cdot)\|_{L^q}\\
+g_{1}(\tau)^{-1}\|v(\tau,\cdot)\|_{L^q} &+ g_{1,\sigma_2}(\tau)^{-1}\big\||D|^{\sigma_2} v(\tau,\cdot)\big\|_{L^q}+ g_{1,2\sigma_2}(\tau)^{-1}\big\||D|^{2\sigma_2} v(\tau,\cdot)\big\|_{L^q}\\
&+ g_{2}(\tau)^{-1}\|v_t(\tau,\cdot)\|_{L^q} \Big),
\end{align*}
where
\begin{align}
&f_{1}(\tau)= (1+\tau)^{-\frac{n}{2\sigma_1}(1-\frac{1}{r})+ \e(p_1,\sigma_2)+ \frac{\kappa_1}{2}},\,\,\, f_{1,2\sigma_1}(\tau)= (1+\tau)^{-\frac{n}{2\sigma_1}(1-\frac{1}{r})-1+ \e(p_1,\sigma_2) + \frac{\kappa_1}{2}}, \label{pt4.11}\\
&f_{1,\sigma_1}(\tau)=f_{2}(\tau)=(1+\tau)^{-\frac{n}{2\sigma_1}(1-\frac{1}{r})- \frac{1}{2} + \e(p_1,\sigma_2) + \frac{\kappa_1}{2}}, \label{pt4.12} \\
&g_{1}(\tau)= (1+\tau)^{-\frac{n}{2\sigma_2}(1-\frac{1}{r})+ \frac{\kappa_1}{2}},\,\,\, g_{1,2\sigma_2}(\tau)= (1+\tau)^{-\frac{n}{2\sigma_2}(1-\frac{1}{r})-1+ \frac{\kappa_1}{2}}, \label{pt4.21}\\
&g_{1,\sigma_2}(\tau)=g_{2}(\tau)=(1+\tau)^{-\frac{n}{2\sigma_2}(1-\frac{1}{r})- \frac{1}{2}+ \frac{\kappa_1}{2}}, \label{pt4.22}
\end{align}
For all $t>0$ we define the following operator $N: \quad (u,v) \in X(t) \longrightarrow N(u,v) \in X(t)$
$$N(u,v)(t,x)= \big(u^{ln}(t,x)+ u^{nl}(t,x), v^{ln}(t,x)+ v^{nl}(t,x)\big). $$
If we prove the operator $N$ satisfying the following two inequalities:
\begin{align}
&\|N(u,v)\|_{X(t)} \lesssim \|(u_0,u_1)\|_{\mathcal{A}^{1,s_1}_{m,q}}+ \|(v_0,v_1)\|_{\mathcal{A}^{2,s_2}_{m,q}}+ \|(u,v)\|^{p_1}_{X(t)}+ \|(u,v)\|^{p_2}_{X(t)}, \label{pt4.3}\\
&\|N(u,v)-N(\bar{u},\bar{v})\|_{X(t)} \nonumber\\
&\qquad \qquad \qquad \lesssim \|(u,v)-(\bar{u},\bar{v})\|_{X(t)} \Big(\|(u,v)\|^{p_1-1}_{X(t)}+ \|(\bar{u},\bar{v})\|^{p_1-1}_{X(t)}+ \|(u,v)\|^{p_2-1}_{X(t)}+ \|(\bar{u},\bar{v})\|^{p_2-1}_{X(t)}\Big), \label{pt4.4}
\end{align}
then we may conclude local (in time) existence results of large data solutions and global (in time) existence results of small data solutions as well by applying Banach's fixed point theorem. \medskip

\noindent In the first step, from the definition of the norm in $X(t)$ we plug $a= \sigma_1,\, \sigma_2,\, 2\sigma_1,\, 2\sigma_2$ into the statements from Corollary \ref{hq2.1} to derive 
\begin{equation}
\big\|(u^{ln}, v^{ln})\big\|_{X(t)} \lesssim \|(u_0,u_1)\|_{\mathcal{A}^{1}_{m,q}}+ \|(v_0,v_1)\|_{\mathcal{A}^{2}_{m,q}}. \label{pt4.5}
\end{equation}
Therefore, in order to prove the proof of (\ref{pt4.3}) it is suitable to indicate the following inequality:
\begin{equation}
\big\|(u^{nl}, v^{nl})\big\|_{X(t)} \lesssim \|(u,v)\|^{p_1}_{X(t)}+ \|(u,v)\|^{p_2}_{X(t)}. \label{pt4.31}
\end{equation}
Now let us prove the inequality (\ref{pt4.31}). To estimate for $u^{nl}$ and some its derivatives, we use the $(L^m \cap L^q)- L^q$ estimates if $\tau \in [0,t/2]$ and the $L^q-L^q$ estimates if $\tau \in [t/2,t]$ from Corollary \ref{hq2.1} to get the following estimates for $k=0,1,2$:
\begin{align*}
\big\||D|^{k\sigma_1} u^{nl}(t,\cdot)\big\|_{L^q} &\lesssim \int_0^{t/2}(1+t-\tau)^{\frac{\gamma+1}{2} -\frac{n}{2\sigma_1}(1-\frac{1}{r})- \frac{k}{2}}\big\||v(\tau,\cdot)|^{p_1}\big\|_{L^m \cap L^q}d\tau\\
&\qquad + \int_{t/2}^t (1+t-\tau)^{\frac{1}{2}(3+[\frac{n}{2}])- \frac{k}{2}}\big\||v(\tau,\cdot)|^{p_1}\big\|_{L^q}d\tau.
\end{align*}
Hence, it is reasonable to control $|v(\tau,x)|^{p_1}$ in $L^m \cap L^q$ and $L^q$. We derive
$$\big\||v(\tau,\cdot)|^{p_1}\big\|_{L^m \cap L^q} \lesssim \|v(\tau,\cdot)\|^{p_1}_{L^{mp_1}}+ \|v(\tau,\cdot)\|^{p_1}_{L^{qp_1}},\, \text{ and }\big\||v(\tau,\cdot)|^{p_1}\big\|_{L^q}= \|v(\tau,\cdot)\|^{p_1}_{L^{qp_1}}.$$
Applying the fractional Gagliardo-Nirenberg inequality leads to
\begin{align*}
\big\||v(\tau,\cdot)|^{p_1}\big\|_{L^m \cap L^q} &\lesssim (1+\tau)^{p_1\big(-\frac{n}{2\sigma_2}(\frac{1}{m}-\frac{1}{mp_1})+ \frac{\kappa_1}{2}\big)}\|(u,v)\|^{p_1}_{X(\tau)}, \\
\big\||v(\tau,\cdot)|^{p_1}\big\|_{L^q} &\lesssim (1+\tau)^{p_1\big(-\frac{n}{2\sigma_2}(\frac{1}{m}-\frac{1}{qp_1})+ \frac{\kappa_1}{2}\big)}\|(u,v)\|^{p_1}_{X(\tau)}.
\end{align*}
Here the conditions (\ref{GN11A1}) to (\ref{GN11A3}) are fulfilled for $p_1$. Consequently, we obtain
\begin{align*}
\big\||D|^{k\sigma_1} u^{nl}(t,\cdot)\big\|_{L^q} &\lesssim (1+t)^{\frac{\gamma+1}{2} -\frac{n}{2\sigma_1}(1-\frac{1}{r})- \frac{k}{2}}\|(u,v)\|^{p_1}_{X(t)} \int_0^{t/2}(1+\tau)^{p_1\big(-\frac{n}{2\sigma_2}(\frac{1}{m}-\frac{1}{mp_1})+ \frac{\kappa_1}{2}\big)} d\tau\\
&\qquad \quad + (1+t)^{p_1\big(-\frac{n}{2\sigma_2}(\frac{1}{m}-\frac{1}{qp_1})+ \frac{\kappa_1}{2}\big)}\|(u,v)\|^{p_1}_{X(t)} \int_{t/2}^t (1+t-\tau)^{\frac{1}{2}(3+[\frac{n}{2}])- \frac{k}{2}}d\tau,
\end{align*}
where we used $(1+t-\tau) \approx (1+t) \text{ for any }\tau \in [0,t/2] \text{ and } (1+\tau) \approx (1+t) \text{ for any }\tau \in [t/2,t] $. Since the condition $p_1 \le 1+\frac{2m\sigma_2(1+\kappa_1)}{n- 2m\sigma_2\kappa_1}$ holds, the term $(1+\tau)^{p_1\big(-\frac{n}{2\sigma_2}(\frac{1}{m}-\frac{1}{mp_1})+ \frac{\kappa_1}{2}\big)}$ is not integrable. As a result, we get
\begin{align*}
&(1+t)^{\frac{\gamma+1}{2} -\frac{n}{2\sigma_1}(1-\frac{1}{r})- \frac{k}{2}} \int_0^{t/2}(1+\tau)^{p_1\big(-\frac{n}{2\sigma_2}(\frac{1}{m}-\frac{1}{mp_1})+ \frac{\kappa_1}{2}\big)} d\tau\\ 
&\qquad \lesssim (1+t)^{-\frac{n}{2\sigma_1}(1-\frac{1}{r})- \frac{k}{2}+\e(p_1,\sigma_2)+ \frac{\kappa_1}{2}}.
\end{align*}
In addition, we also note that $\frac{1}{2}(3+[\frac{n}{2}])- \frac{k}{2}> 0$ for $k=0,1,2$. Therefore, we arrive at
\begin{align*}
&(1+t)^{p_1\big(-\frac{n}{2\sigma_2}(\frac{1}{m}-\frac{1}{qp_1})+ \frac{\kappa_1}{2}\big)} \int_{t/2}^t (1+t-\tau)^{\frac{1}{2}(3+[\frac{n}{2}])- \frac{k}{2}}d\tau\\ 
&\qquad \lesssim (1+t)^{-\frac{n}{2\sigma_1}(1-\frac{1}{r})- \frac{k}{2}+\e(p_1,\sigma_2)+ \frac{\kappa_1}{2}},
\end{align*}
since $\sigma_1 \ge \sigma_2$. From both the above estimates, we may conclude the following estimate for $k=0,1,2$:
$$\big\||D|^{k\sigma_1} u^{nl}(t,\cdot)\big\|_{L^q} \lesssim (1+t)^{-\frac{n}{2\sigma_1}(1-\frac{1}{r})- \frac{k}{2}+\e(p_1,\sigma_2)+ \frac{\kappa_1}{2}} \|(u,v)\|^{p_1}_{X(t)}. $$
In the analogous way we also derive
$$\big\|u_t^{nl}(t,\cdot)\big\|_{L^q} \lesssim (1+t)^{-\frac{n}{2\sigma_1}(1-\frac{1}{r})- \frac{1}{2}+\e(p_1,\sigma_2)+ \frac{\kappa_1}{2}} \|(u,v)\|^{p_1}_{X(t)}. $$
Similarly, we obtain the following estimate for $k=0,1,2$:
$$\big\||D|^{k\sigma_2} v^{nl}(t,\cdot)\big\|_{L^2} \lesssim (1+t)^{-\frac{n}{2\sigma_2}(1-\frac{1}{r})- \frac{k}{2}+ \frac{\kappa_1}{2}}\|(u,v)\|^{p_2}_{X(t)}, $$
and
$$\big\|v_t^{nl}(t,\cdot)\big\|_{L^2} \lesssim (1+t)^{-\frac{n}{2\sigma_2}(1-\frac{1}{r})- \frac{1}{2} + \frac{\kappa_1}{2}}\|(u,v)\|^{p_2}_{X(t)}, $$
where the conditions (\ref{GN11A1}) to (\ref{exponent11A}) are satisfied for $p_2$. From the definition of the norm in $X(t)$, we may conclude immediately the inequality (\ref{pt4.31}). \medskip

\noindent In the second step, let us prove the estimate (\ref{pt4.4}). We obtain for two elements $(u,v)$ and $(\bar{u},\bar{v})$ from $X(t)$ as follows:
$$N(u,v)(t,x)- N(\bar{u},\bar{v})(t,x)= \big(u^{nl}(t,x)- \bar{u}^{nl}(t,x), v^{nl}(t,x)- \bar{v}^{nl}(t,x)\big). $$
We use again the $(L^m \cap L^q)- L^q$ estimates if $\tau \in [0,t/2]$ and the $L^q-L^q$ estimates if $\tau \in [t/2,t]$ from Corollary \ref{hq2.1} to get the following estimate for $k=0,1,2$:
\begin{align*}
\big\||D|^{k\sigma_1}\big(u^{nl}- \bar{u}^{nl}\big)(t,\cdot)\big\|_{L^q} &\lesssim \int_0^{t/2}(1+t-\tau)^{\frac{\gamma+1}{2} -\frac{n}{2\sigma_1}(1-\frac{1}{r})- \frac{k}{2}}\big\||v(\tau,\cdot)|^{p_1}- \bar{v}(\tau,\cdot)|^{p_1}\big\|_{L^m \cap L^q}d\tau \\
&\qquad + \int_{t/2}^t (1+t-\tau)^{\frac{1}{2}(3+[\frac{n}{2}])- \frac{k}{2}}\big\||v(\tau,\cdot)|^{p_1}- \bar{v}(\tau,\cdot)|^{p_1}\big\|_{L^q}d\tau, \\
\big\|\big(u_t^{nl}- \bar{u}_t^{nl}\big)(t,\cdot)\big\|_{L^q} &\lesssim \int_0^{t/2}(1+t-\tau)^{\frac{\gamma}{2} -\frac{n}{2\sigma_1}(1-\frac{1}{r})}\big\||v(\tau,\cdot)|^{p_1}- \bar{v}(\tau,\cdot)|^{p_1}\big\|_{L^m \cap L^q}d\tau \\
&\qquad + \int_{t/2}^t (1+t-\tau)^{\frac{1}{2}(2+[\frac{n}{2}])}\big\||v(\tau,\cdot)|^{p_1}- \bar{v}(\tau,\cdot)|^{p_1}\big\|_{L^q}d\tau,
\end{align*}
and
\begin{align*}
\big\||D|^{k\sigma_2}\big(v^{nl}- \bar{v}^{nl}\big)(t,\cdot)\big\|_{L^q} &\lesssim \int_0^{t/2}(1+t-\tau)^{\frac{\gamma+1}{2} -\frac{n}{2\sigma_2}(1-\frac{1}{r})- \frac{k}{2}}\big\||u(\tau,\cdot)|^{p_2}- \bar{u}(\tau,\cdot)|^{p_2}\big\|_{L^m \cap L^q}d\tau \\
&\qquad + \int_{t/2}^t (1+t-\tau)^{\frac{1}{2}(3+[\frac{n}{2}])- \frac{k}{2}}\big\||u(\tau,\cdot)|^{p_2}- \bar{u}(\tau,\cdot)|^{p_2}\big\|_{L^q}d\tau, \\
\big\|\big(v_t^{nl}- \bar{v}_t^{nl}\big)(t,\cdot)\big\|_{L^q} &\lesssim \int_0^{t/2}(1+t-\tau)^{\frac{\gamma}{2} -\frac{n}{2\sigma_2}(1-\frac{1}{r})}\big\||u(\tau,\cdot)|^{p_2}- \bar{u}(\tau,\cdot)|^{p_2}\big\|_{L^m \cap L^q}d\tau \\
&\qquad + \int_{t/2}^t (1+t-\tau)^{\frac{1}{2}(2+[\frac{n}{2}])}\big\||u(\tau,\cdot)|^{p_2}- \bar{u}(\tau,\cdot)|^{p_2}\big\|_{L^q}d\tau.
\end{align*}
By using H\"{o}lder's inequality we arrive at
\begin{align*}
\big\||v(\tau,\cdot)|^{p_1}- |\bar{v}(\tau,\cdot)|^{p_1}\big\|_{L^q}& \lesssim \|v(\tau,\cdot)- \bar{v}(\tau,\cdot)\|_{L^{qp_1}} \big(\|v(\tau,\cdot)\|^{p_1-1}_{L^{qp_1}}+ \|\bar{v}(\tau,\cdot)\|^{p_1-1}_{L^{qp_1}}\big),\\
\big\||v(\tau,\cdot)|^{p_1}- |\bar{v}(\tau,\cdot)|^{p_1}\big\|_{L^m}& \lesssim \|v(\tau,\cdot)- \bar{v}(\tau,\cdot)\|_{L^{mp_1}} \big(\|v(\tau,\cdot)\|^{p_1-1}_{L^{mp_1}}+ \|\bar{v}(\tau,\cdot)\|^{p_1-1}_{L^{mp_1}}\big),\\
\big\||u(\tau,\cdot)|^{p_2}- |\bar{u}(\tau,\cdot)|^{p_2}\big\|_{L^q}& \lesssim \|u(\tau,\cdot)- \bar{u}(\tau,\cdot)\|_{L^{qp_2}} \big(\|u(\tau,\cdot)\|^{p_2-1}_{L^{qp_2}}+ \|\bar{u}(\tau,\cdot)\|^{p_2-1}_{L^{qp_2}}\big),\\
\big\||u(\tau,\cdot)|^{p_2}- |\bar{u}(\tau,\cdot)|^{p_2}\big\|_{L^m}& \lesssim \|u(\tau,\cdot)- \bar{u}(\tau,\cdot)\|_{L^{mp_2}} \big(\|u(\tau,\cdot)\|^{p_2-1}_{L^{mp_2}}+ \|\bar{u}(\tau,\cdot)\|^{p_2-1}_{L^{mp_2}}\big).
\end{align*}
In the same way as we proved (\ref{pt4.3}), employing the fractional Gagliardo-Nirenberg inequality to the terms
$$ \|v(\tau,\cdot)-\bar{v}(\tau,\cdot)\|_{L^{\eta_1}},\,\, \|u(\tau,\cdot)-\bar{u}(\tau,\cdot)\|_{L^{\eta_2}},\,\, \|v(\tau,\cdot)\|_{L^{\eta_1}},\,\, \|\bar{v}(\tau,\cdot)\|_{L^{\eta_1}},\,\, \|u(\tau,\cdot)\|_{L^{\eta_2}},\,\, \|\bar{u}(\tau,\cdot)\|_{L^{\eta_2}} $$
with $\eta_1=qp_1$ or $\eta_1=mp_1$, and $\eta_2=qp_2$ or $\eta_2=mp_2$ we may conclude the inequality (\ref{pt4.4}). Summarizing, the proof of Theorem \ref{dl1.1} is completed.

\subsection{Proof of Theorem \ref{dl1.2}}
We follow the approach of Theorem \ref{dl1.1} with minor modifications in the steps of our proof. We also introduce both spaces for the data and the solutions as in Theorem \ref{dl1.1}, where the weights (\ref{pt4.11}) to (\ref{pt4.22}) are modified in the following way: 
\begin{align*}
&f_{1}(\tau)= (1+\tau)^{-\frac{n}{2\sigma_1}(1-\frac{1}{r})+ \frac{\kappa_2}{2}},\,\,\, f_{1,2\sigma_1}(\tau)= (1+\tau)^{-\frac{n}{2\sigma_1}(1-\frac{1}{r})-1+ \frac{\kappa_2}{2}}, \\
&f_{1,\sigma_1}(\tau)=f_{2}(\tau)=(1+\tau)^{-\frac{n}{2\sigma_1}(1-\frac{1}{r})- \frac{1}{2} + \frac{\kappa_2}{2}}, \\
&g_{1}(\tau)= (1+\tau)^{-\frac{n}{2\sigma_2}(1-\frac{1}{r})+ \e(p_2,\sigma_1)+ \frac{\kappa_2}{2}},\,\,\, g_{1,2\sigma_2}(\tau)= (1+\tau)^{-\frac{n}{2\sigma_2}(1-\frac{1}{r})-1+ \e(p_2,\sigma_1)+ \frac{\kappa_2}{2}}, \\
&g_{1,\sigma_2}(\tau)=g_{2}(\tau)=(1+\tau)^{-\frac{n}{2\sigma_2}(1-\frac{1}{r})- \frac{1}{2} + \e(p_2,\sigma_1)+ \frac{\kappa_2}{2}},
\end{align*}
Then, repeating some steps of the proofs we did in Theorem \ref{dl1.1} we may complete the proof of Theorem \ref{dl1.2}.

\subsection{Proof of Theorem \ref{dl1.1B}}
The proof of Theorem \ref{dl1.1B} is similar to the proof of Theorem \ref{dl1.1}. We also introduce both spaces for the data and the solutions as in Theorem \ref{dl1.1}, where the weights (\ref{pt4.11}) to (\ref{pt4.22}) are modified in the following way:
\begin{align*}
&f_{1}(\tau)= (1+\tau)^{-\frac{n}{2\sigma_1}(1-\frac{1}{r})+ \frac{\gamma+1}{2}},\,\,\, f_{1,2\sigma_1}(\tau)= (1+\tau)^{-\frac{n}{2\sigma_1}(1-\frac{1}{r})+ \frac{\gamma-1}{2}}, \\
&f_{1,\sigma_1}(\tau)=f_{2}(\tau)=(1+\tau)^{-\frac{n}{2\sigma_1}(1-\frac{1}{r})+ \frac{\gamma}{2}}, \\
&g_{1}(\tau)= (1+\tau)^{-\frac{n}{2\sigma_2}(1-\frac{1}{r})+ \frac{\gamma+1}{2}},\,\,\, g_{1,2\sigma_2}(\tau)= (1+\tau)^{-\frac{n}{2\sigma_2}(1-\frac{1}{r})+ \frac{\gamma-1}{2}}, \\
&g_{1,\sigma_2}(\tau)=g_{2}(\tau)=(1+\tau)^{-\frac{n}{2\sigma_2}(1-\frac{1}{r})+ \frac{\gamma}{2}},
\end{align*}
Here we notice that due to the condition (\ref{exponent11B}), the terms
$$(1+\tau)^{p_1\big(-\frac{n}{2\sigma_2}(\frac{1}{m}-\frac{1}{mp_1})+ \frac{\kappa_1}{2}\big)} \text{ and } (1+\tau)^{p_2\big(-\frac{n}{2\sigma_2}(\frac{1}{m}-\frac{1}{mp_2})+ \frac{\kappa_1}{2}\big)}$$
are integrable. Then, repeating some of the arguments as we did in the proof of Theorem \ref{dl1.1} we may complete the proof of Theorem \ref{dl1.1B}.

\subsection{Proof of Theorem \ref{dl1.2B}}
The proof of Theorem \ref{dl1.2B} is similar to the proof of Theorem \ref{dl1.1}. We also introduce both spaces for the data and the solutions as in Theorem \ref{dl1.1}, where the weights are the same as those in Theorem \ref{dl1.1B}. Here we notice that due to the condition (\ref{exponent12B}), the terms
$$(1+\tau)^{p_1\big(-\frac{n}{2\sigma_1}(\frac{1}{m}-\frac{1}{mp_1})+ \frac{\kappa_2}{2}\big)} \text{ and } (1+\tau)^{p_2\big(-\frac{n}{2\sigma_1}(\frac{1}{m}-\frac{1}{mp_2})+ \frac{\kappa_2}{2}\big)}$$
are integrable. Then, repeating some of the arguments as we did in the proof of Theorem \ref{dl1.1} we may complete the proof of Theorem \ref{dl1.2B}.\medskip

\noindent \textbf{Acknowledgments}\medskip

\noindent The PhD study of MSc. T.A. Dao is supported by Vietnamese Government's Scholarship. The author would like to thank sincerely to Prof. Michael Reissig for valuable discussions and Institute of Applied Analysis for their hospitality. The author is grateful to the referee for his careful reading of the manuscript and for helpful comments.


\end{document}